# Normal and Spherical Curves in Dual Space $D^3$

**Mehmet Önder[1], H. Hüseyin Uğurlu[2]**

[1]Celal Bayar University, Faculty of Science and Arts, Mathematics Department, Muradiye Campus, Muradiye, Manisa, Turkey. E-mail: mehmet.onder@cbu.edu.tr

[2]Gazi University, Faculty of Education, Department of Secondary Education Science and Mathematics Teaching, Mathematics Teaching Programme, Ankara, Turkey. E-mail: hugurlu@gazi.edu.tr

**Abstract:** In this paper, we give definitions and characterizations of normal and spherical curves in the dual space $D^3$. We show that normal curves are also spherical curves in $D^3$.

## 1. Introduction

As it is well known, the set of orientable lines in $E^3$ admits a bijection with the set of vectors $(\vec{a},\vec{a}^*)$ in $\{(\vec{a},\vec{a}^*) \in E^3 \times E^3 : \langle \vec{a},\vec{a} \rangle = 1, \ \langle \vec{a},\vec{a}^* \rangle = 0\}$. This bijection is obtained as follows; such an oriented line $L$ is determined by a point $p \in L$ and a unit vector $\vec{a}$. Then, one can define $\vec{a}^* = \vec{p} \times \vec{a}$. The value of $\vec{a}^*$ does not depend on the point $p$, because any other point $q$ in $L$ is $\vec{q} = \vec{p} + \lambda \vec{a}$ and then $\vec{a}^* = \vec{p} \times \vec{a} = \vec{q} \times \vec{a}$. Reciprocally, given such a pair $(\vec{a},\vec{a}^*)$ one recovers the line $L$ as $L = \{(\vec{a} \times \vec{a}^*) + \lambda \vec{a} : \vec{a}, \vec{a}^* \in E^3, \lambda \in IR\}$, written in parametric equations. The above mentioned set of vectors is called the dual unit sphere $\tilde{S}^2$, because it has a relation with dual numbers. Since each dual unit vector corresponds to a line of $E^3$, there is a one-to-one correspondence between the points of the dual unit sphere $\tilde{S}^2$ and the oriented lines of $E^3$. This correspondence is known as E. Study Mapping [2,3,10]. As a consequence of that, a differentiable curve lying fully on dual unit sphere in dual space $D^3$ represents a ruled surface which is a surface generated by moving of a line $L$ along a curve $\vec{\alpha}(s)$ in $E^3$ and has the parametrization $\vec{r}(s,u) = \vec{\alpha}(s) + u\vec{l}(s)$, where $\vec{\alpha}(s)$ is called indicatrix curve and $\vec{l}(s)$ is the direction vector of the generating line which passes by the point $\vec{\alpha}(s)$. These straight lines are called rulings or rectilinear generators of the surface.

In the Euclidean space $E^3$, the curve $\alpha : I \subset IR \to E^3$ for which the position vector $\vec{\alpha}(s)$ always lies in its rectifying plane, is for simplicity called *rectifying curve*. Similarly, the curve $\vec{\alpha}(s)$ for which the position vector $\vec{\alpha}(s)$ always lies in its osculating plane (respectively normal plane) is called *osculating curve (respectively normal curve)* [1]. The characterizations of rectifying and normal curves in Euclidean space and Minkowski space have been given some authors [1,4,5,6]. The corresponding characterizations of the rectifying dual curves in dual space have been given by Yücesan, Ayyıldız and Çöken [9]. Furthermore, Önder has defined and studied dual timelike normal and dual timelike spherical curves in dual Minkowski space $D_1^3$ [7].

In this paper, we give the characterizations of normal curves and spherical curves in dual space $D^3$. Furthermore, we show that every dual normal curve is also a dual spherical curve in dual space $D^3$.

## 2. Preliminaries

Dual numbers were introduced by W.K. Clifford (1845-1879). A dual number has the form $\bar{a} = (a, a^*) = a + \varepsilon a^*$ where $a$ and $a^*$ are real numbers and $\varepsilon = (0,1)$ is called dual unit. We denote the set of dual numbers by $D$, i.e.,

$$D = \{\bar{a} = a + \varepsilon a^* : a, a^* \in \mathbb{R}, \varepsilon^2 = 0\}.$$

In $D$ two inner operations are defined as follows:

**i)** An operation "+" defined as
$$\bar{a} + \bar{b} = (a, a^*) + (b, b^*) = (a+b,\ a^* + b^*) = a + b + \varepsilon(a^* + b^*),$$
and is called the addition in $D$.

**ii)** An operation "·" defined by
$$\bar{a} \cdot \bar{b} = \bar{a}\bar{b} = (a, a^*)(b, b^*) = (ab,\ ab^* + a^*b) = ab + \varepsilon(ab^* + a^*b),$$
and is called the multiplication in $D$.

If $\bar{a} = a + \varepsilon a^*$ is a dual number, the real number $a$ is called the real part of $\bar{a}$ and the real number $a^*$ is called the dual part of $\bar{a}$. For definition of the multiplication $\varepsilon^2 = 0$.

The set $D$ of dual numbers with the above operations is a commutative ring.

The dual number $\bar{a} = a + \varepsilon a^*$ divide by dual number $\bar{b} = b + \varepsilon b^*$, with $b \neq 0$, is defined by

$$\frac{\bar{a}}{\bar{b}} = \frac{a}{b} + \varepsilon \frac{a^*b - ab^*}{b^2}.$$

Now let $f$ be a differentiable function with dual variable $\bar{x} = x + \varepsilon x^*$. Then the Maclaurin series generated by $f$ is given by

$$f(\bar{x}) = f(x + \varepsilon x^*) = f(x) + \varepsilon x^* f'(x), \tag{1}$$

where $f'(x)$ is derivative of $f(x)$ with respect to $x$. Then we have

$$\sin(x + \varepsilon x^*) = \sin(x) + \varepsilon x^* \cos(x),$$
$$\cos(x + \varepsilon x^*) = \cos(x) - \varepsilon x^* \sin(x).$$

Let $D^3$ be the set of all triples of dual numbers, i.e.,
$$D^3 = \{\tilde{a} = (\bar{a}_1, \bar{a}_2, \bar{a}_3) : \bar{a}_i \in D, i = 1, 2, 3\},$$

which is called dual space, and it is a module over the ring $D$. The elements of $D^3$ are called dual vectors. A dual vector $\tilde{a}$ may be expressed in the form $\tilde{a} = \vec{a} + \varepsilon \vec{a}^* = (\vec{a}, \vec{a}^*)$, where $\vec{a} = (a_1, a_2, a_3)$ and $\vec{a}^* = (a_1^*, a_2^*, a_3^*)$ are vectors of $\mathbb{R}^3$. Then for any dual vectors $\tilde{a} = \vec{a} + \varepsilon \vec{a}^*$ and $\tilde{b} = \vec{b} + \varepsilon \vec{b}^*$ in $D^3$, the scalar product and the vector product are defined by

$$g(\tilde{a}, \tilde{b}) = g(\vec{a}, \vec{b}) + \varepsilon\left(g(\vec{a}, \vec{b}^*) + g(\vec{a}^*, \vec{b})\right),$$

and

$$\tilde{a} \times \tilde{b} = \vec{a} \times \vec{b} + \varepsilon\left(\vec{a} \times \vec{b}^* + \vec{a}^* \times \vec{b}\right),$$

respectively, where $g(\vec{a}, \vec{b})$ is the inner product of the vectors $\vec{a}$ and $\vec{a}^*$ in $\mathbb{R}^3$.

The norm of a dual vector $\tilde{a}$ is given by

$$\|\tilde{a}\| = \sqrt{g(\tilde{a},\tilde{a})} = \|\vec{a}\| + \varepsilon \frac{g(\vec{a},\vec{a}^*)}{\|\vec{a}\|}.$$

A dual vector $\tilde{a}$ with norm 1 is called dual unit vector. The set of dual unit vectors is denoted by

$$\tilde{S}^2 = \{\tilde{a} = (\overline{a}_1, \overline{a}_2, \overline{a}_3) \in D^3 : g(\tilde{a},\tilde{a}) = (1,0)\},$$

and called dual unit sphere. Similarly, the set of arbitrary dual vectors

$$\tilde{S}^2(\overline{r}) = \{\tilde{a} = (\overline{a}_1, \overline{a}_2, \overline{a}_3) \in D^3 : g(\tilde{a}-\overline{c}, \tilde{a}-\overline{c}) = \overline{r}^2, \overline{c} \in D\},$$

is a dual sphere with radius $\overline{r}$ and center $\overline{c}$ (For details [2,3,10]).

E. Study used dual numbers and dual vectors in his research on the geometry of lines and kinematics. He devoted special attention to the representation of directed lines by dual unit vectors and defined the mapping that is known by his name: There exists a one-to-one correspondence between the vectors of dual unit sphere $\tilde{S}^2$ and the directed lines of space of lines $\mathbb{R}^3$ [2,10]. By the aid of this correspondence, the properties of the spatial motion of a line can be derived. Hence, the geometry of the ruled surface is represented by the geometry of dual curve lying fully on the dual unit sphere in $D^3$.

The angle $\overline{\theta} = \theta + \varepsilon\theta^*$ between two dual unit vectors $\tilde{a}, \tilde{b}$ is called dual angle and defined by

$$g(\tilde{a},\tilde{b}) = \cos\overline{\theta} = \cos\theta - \varepsilon\theta^* \sin\theta.$$

By considering the E. Study Mapping, the geometric interpretation of dual angle is that $\theta$ is the real angle between lines $L_1$, $L_2$ corresponding to the dual unit vectors $\tilde{a}, \tilde{b}$ respectively, and $\theta^*$ is the shortest distance between those lines.

Let $\vec{\alpha}(t) = (\alpha_1(t), \alpha_2(t), \alpha_3(t),)$ and $\vec{\alpha}^*(t) = (\alpha_1^*(t), \alpha_2^*(t), \alpha_3^*(t))$ be real valued curves in the space $E^3$. Then $\tilde{\alpha}(t) = \vec{\alpha}(t) + \varepsilon\vec{\alpha}^*(t)$ is a curve in dual space $D^3$ and is called *dual space curve* or *curve in dual space*. If the real valued functions $\alpha_i(t)$ and $\alpha_i^*(t)$, $(1 \leq i \leq 3)$ are differentiable then the dual space curve

$$\tilde{\alpha}: I \subset \mathrm{IR} \to D^3$$
$$t \to \tilde{\alpha}(t) = (\alpha_1(t) + \varepsilon\alpha_1^*(t), \alpha_2(t) + \varepsilon\alpha_2^*(t), \alpha_3(t) + \varepsilon\alpha_3^*(t)) \quad (2)$$
$$= \vec{\alpha}(t) + \varepsilon\vec{\alpha}^*(t)$$

is differentiable in dual space $D^3$. The real part $\vec{\alpha}(t)$ of the dual space curve $\tilde{\alpha} = \tilde{\alpha}(t)$ is called *indicatrix*. The dual arc-length of the dual space curve $\tilde{\alpha}(t)$ from $t_1$ to $t$ is defined by

$$\overline{s} = \int_{t_1}^{t} \|\tilde{\alpha}'(t)\| dt = \int_{t_1}^{t} \|\vec{\alpha}'(t)\| dt + \varepsilon \int_{t_1}^{t} g\left(\vec{T}, \vec{\alpha}^{*'}(t)\right) dt = s + \varepsilon s^* \quad (3)$$

where $\vec{T}$ is unit tangent vector of the indicatrix $\vec{\alpha}(t)$ which is a real space curve in $E^3$ [7,9]. From the last equality we can write

$$\frac{d\overline{s}}{ds} = 1 + \varepsilon\Delta \quad (4)$$

where $\Delta = g\left(\vec{T}, \vec{\alpha}^{*'}\right)$.

Now, let $\tilde{\alpha}(s)$ be a $C^4$-dual space curve with arc length parameter $s$ of the indicatrix. Then dual unit tangent of $\tilde{\alpha}$ is defined by

$$\frac{d\tilde{\alpha}}{d\bar{s}} = \tilde{T}.$$

Since $\tilde{T}$ is unit, from the differentiation of $\tilde{T}$ with respect to dual arc length parameter $\bar{s}$ we have

$$\tilde{T}' = \frac{d\tilde{T}}{d\bar{s}} = \frac{d^2\tilde{\alpha}}{d\bar{s}^2} = \bar{k}_1 \tilde{N},$$

where $\bar{k}_1 = \bar{k}_1(s)$ is called dual curvature which is the norm of $\tilde{T}'$. We impose the restriction that $\bar{k}_1(s)$ is never pure dual. Then the dual unit vector $\tilde{N} = (1/\bar{k}_1)\tilde{T}'$ is called dual unit principal normal vector of $\tilde{\alpha}$. The dual unit vector $\tilde{B}$ defined by $\tilde{B} = \tilde{T} \times \tilde{N}$ is called dual unit binormal vector of $\tilde{\alpha}$. Then the dual frame denoted by $\{\tilde{T}(\bar{s}), \tilde{N}(\bar{s}), \tilde{B}(\bar{s})\}$ is called moving dual Frenet frame along the dual space curve $\tilde{\alpha}(s)$ in the dual space $D^3$. Then for the curve $\tilde{\alpha}$ the dual Frenet formulae are given by

$$\frac{d}{d\bar{s}}\begin{bmatrix}\tilde{T}\\ \tilde{N}\\ \tilde{B}\end{bmatrix} = \begin{bmatrix}0 & \bar{k}_1 & 0\\ -\bar{k}_1 & 0 & \bar{k}_2\\ 0 & -\bar{k}_2 & 0\end{bmatrix}\begin{bmatrix}\tilde{T}\\ \tilde{N}\\ \tilde{B}\end{bmatrix}, \tag{5}$$

where $\bar{k}_2 = \bar{k}_2(s)$ is called the dual torsion of $\tilde{\alpha}$. Similar to the Euclidean case, the planes spanned by $\{\tilde{T}, \tilde{N}\}$, $\{\tilde{T}, \tilde{B}\}$ and $\{\tilde{N}, \tilde{B}\}$ are called the dual osculating plane, the dual rectifying plane and the dual normal plane, respectively [9].

### 3. Dual Normal Curves in $D^3$

In this section, we give definition and characterizations of normal curves in dual space $D^3$.

**Definition 3.1.** Let $\tilde{\alpha} = \tilde{\alpha}(s)$ be a curve in dual space $D^3$. $\tilde{\alpha}$ is called *dual normal curve (or normal curve in $D^3$)* if

$$\tilde{\alpha}(s) = \bar{\lambda}(s)\tilde{N}(s) + \bar{\mu}(s)\tilde{B}(s),$$

where $\bar{\lambda}(s)$ and $\bar{\mu}(s)$ are dual differentiable functions of arc-length parameter $s$ of indicatrix.

Let $\tilde{\alpha} = \tilde{\alpha}(s)$ be a dual curve with dual curvature $\bar{k}_1(s)$ and dual torsion $\bar{k}_2(s)$, and suppose that $\bar{k}_2$ also is nowhere pure dual. Then we have the followings:

**Theorem 3.1.** Let $\tilde{\alpha} = \tilde{\alpha}(s)$ be a normal curve in $D^3$ with dual curvatures $\bar{k}_1(s) \neq 0$, $\bar{k}_2(s) \neq 0$. Then the following statements hold:

*i)* The dual curvatures $\bar{k}_1(s)$ and $\bar{k}_2(s)$ satisfy the following equality

$$\frac{1}{\bar{k}_1(s)} = \bar{c}_1 \cos\left(\int \bar{k}_2(s)d\bar{s}\right) + \bar{c}_2 \sin\left(\int \bar{k}_2(s)d\bar{s}\right).$$

*(ii)* The principal normal and binormal components of position vector of dual space curve $\tilde{\alpha}$ are given by

$$g(\tilde{\alpha}, \tilde{N}) = -\bar{c}_1 \cos\left(\int \bar{k}_2(s)d\bar{s}\right) - \bar{c}_2 \sin\left(\int \bar{k}_2(s)d\bar{s}\right),$$

$$g(\tilde{\alpha}, \tilde{B}) = \bar{c}_1 \sin\left(\int \bar{k}_2(s)d\bar{s}\right) - \bar{c}_2 \cos\left(\int \bar{k}_2(s)d\bar{s}\right),$$

*respectively, where* $\bar{c}_1 = c_1 + \varepsilon c_1^*$, $\bar{c}_2 = c_2 + \varepsilon c_2^* \in D$ *and* $c_1, c_1^*, c_2, c_2^* \in \mathbb{R}$.

*Conversely, if* $\tilde{\alpha} = \tilde{\alpha}(s)$ *is a dual space curve in* $D^3$ *with the curvatures* $\bar{k}_1(s) \neq 0$, $\bar{k}_2(s) \neq 0$ *and one of the statements (i) and (ii) holds, then* $\tilde{\alpha}$ *is a normal curve in* $D^3$ *or congruent to a normal curve.*

**Proof:** Assume that $\tilde{\alpha} = \tilde{\alpha}(s)$ is a dual normal curve in $D^3$. Then, by Definition 3.1, we have

$$\tilde{\alpha}(s) = \bar{\lambda}(s)\tilde{N}(s) + \bar{\mu}(s)\tilde{B}(s),$$

where $\bar{\lambda}(s)$ and $\bar{\mu}(s)$ are dual differentiable functions of arc-length parameter $s$. Differentiating this equality with respect to the dual arc length parameter $\bar{s}$ and using the Frenet equations (5), we obtain

$$-\bar{\lambda}\bar{k}_1 = 1, \quad \bar{\lambda}' - \bar{\mu}\bar{k}_2 = 0, \quad \bar{\lambda}\bar{k}_2 + \bar{\mu}' = 0. \tag{6}$$

From the first and second equations in (6), we get

$$\bar{\lambda} = -\frac{1}{\bar{k}_1}, \quad \bar{\mu} = -\frac{1}{\bar{k}_2}\left(\frac{1}{\bar{k}_1}\right)'. \tag{7}$$

Thus

$$\tilde{\alpha}(s) = -\frac{1}{\bar{k}_1}\tilde{N} - \frac{1}{\bar{k}_2}\left(\frac{1}{\bar{k}_1}\right)'\tilde{B}. \tag{8}$$

Further, from the third equation in (6) and using (7) we find the following differential equation

$$\left[\frac{1}{\bar{k}_2}\left(\frac{1}{\bar{k}_1}\right)'\right]' + \frac{\bar{k}_2}{\bar{k}_1} = 0. \tag{9}$$

Putting $\bar{y}(s) = \dfrac{1}{\bar{k}_1(s)}$ and $\bar{z}(s) = \dfrac{1}{\bar{k}_2(s)}$, equation (9) can be written as

$$\left(\bar{z}(s)\bar{y}'(s)\right)' + \frac{\bar{y}(s)}{\bar{z}(s)} = 0.$$

If we change variables in the above equation as $\bar{t} = \int \dfrac{1}{\bar{z}(s)}d\bar{s} = \int \bar{k}_2(s)d\bar{s}$ then we get

$$\frac{d^2\bar{y}}{d\bar{t}^2} + \bar{y} = 0.$$

The solution of this differential equation is

$$\bar{y} = \bar{c}_1\cos(\bar{t}) + \bar{c}_2\sin(\bar{t}),$$

where $\bar{c}_1, \bar{c}_2 \in D$. Therefore

$$\frac{1}{\bar{k}_1(s)} = \bar{c}_1\cos\left(\int \bar{k}_2(s)d\bar{s}\right) + \bar{c}_2\sin\left(\int \bar{k}_2(s)d\bar{s}\right). \tag{10}$$

Thus the statement (i) is proved. Next, substituting (10) into (7) and (8) we get

$$\begin{cases}\bar{\lambda} = -\bar{c}_1\cos\left(\int \bar{k}_2(s)d\bar{s}\right) - \bar{c}_2\sin\left(\int \bar{k}_2(s)d\bar{s}\right), \\ \bar{\mu} = \bar{c}_1\sin\left(\int \bar{k}_2(s)d\bar{s}\right) - \bar{c}_2\cos\left(\int \bar{k}_2(s)d\bar{s}\right)\end{cases} \tag{11}$$

and

$$\tilde{\alpha} = \left[ -\overline{c}_1 \cos\left(\int \overline{k}_2(s)d\overline{s}\right) - \overline{c}_2 \sin\left(\int \overline{k}_2(s)d\overline{s}\right) \right] \tilde{N}$$
$$+ \left[ \overline{c}_1 \sin\left(\int \overline{k}_2(s)d\overline{s}\right) - \overline{c}_2 \cos\left(\int \overline{k}_2(s)d\overline{s}\right) \right] \tilde{B} \quad (12)$$

respectively. From (12) we find

$$g(\tilde{\alpha}, \tilde{\alpha}) = \overline{c}_1^2 + \overline{c}_2^2, \quad (13)$$

$$g(\tilde{\alpha}, \tilde{N}) = -\overline{c}_1 \cos\left(\int \overline{k}_2(s)d\overline{s}\right) - \overline{c}_2 \sin\left(\int \overline{k}_2(s)d\overline{s}\right), \quad (14)$$

$$g(\tilde{\alpha}, \tilde{B}) = \overline{c}_1 \sin\left(\int \overline{k}_2(s)d\overline{s}\right) - \overline{c}_2 \cos\left(\int \overline{k}_2(s)d\overline{s}\right). \quad (15)$$

Consequently, we have proved (ii).

Conversely, suppose that statement (i) holds. Then we have

$$\frac{1}{\overline{k}_1(s)} = \overline{c}_1 \cos\left(\int \overline{k}_2(s)d\overline{s}\right) + \overline{c}_2 \sin\left(\int \overline{k}_2(s)d\overline{s}\right), \quad \overline{c}_1, \overline{c}_2 \in D.$$

Differentiation of this last equation with respect to $\overline{s}$ gives

$$\left[ \frac{1}{\overline{k}_2} \left( \frac{1}{\overline{k}_1} \right)' \right]' = -\frac{\overline{k}_2}{\overline{k}_1},$$

and by applying dual Frenet formulae (5), we obtain

$$\left[ \tilde{\alpha}(s) - \frac{1}{\overline{k}_1} \tilde{N} - \frac{1}{\overline{k}_2} \left( \frac{1}{\overline{k}_1} \right)' \tilde{B} \right]' = 0.$$

Consequently, $\tilde{\alpha}$ is congruent to a normal curve in $D^3$. Next, assume that statement (ii) holds. Then (13) and (14) are satisfied. Differentiating (13) with respect to $\overline{s}$ and using (14) we find $g(\tilde{\alpha}, \tilde{T}) = 0$, which means that $\tilde{\alpha}$ is dual normal curve, which proves the theorem.

From (10), (14) and (15) we have the following corollaries:

***Corollary 3.1.*** *The real and dual parts of the equation* (10) *are given by*

$$\frac{1}{k_1} = c_1 \cos\left(\int k_2 ds\right) + c_2 \sin\left(\int k_2 ds\right), \quad c_1, c_2 \in \mathbb{R},$$

*and*

$$-\frac{k_1^*}{k_1^2} = \left[ -c_1 \sin\left(\int k_2 ds\right) + c_2 \cos\left(\int k_2 ds\right) \right] \int \left( \Delta k_2 + k_2^* \right) ds$$
$$+ c_1^* \cos\left(\int k_2 ds\right) + c_2^* \sin\left(\int k_2 ds\right).$$

*respectively. Similarly, the real and dual parts of the equations* (14) *and* (15) *are given by*

$$g(\alpha, N) = -c_1 \cos\left(\int k_2 ds\right) - c_2 \sin\left(\int k_2 ds\right)$$

$$g(\alpha, N^*) + g(\alpha^*, N) = \left[ c_1 \sin\left(\int k_2 ds\right) - c_2 \cos\left(\int k_2 ds\right) \right] \int \left( \Delta k_2 + k_2^* \right) ds$$
$$- c_1^* \cos\left(\int k_2 ds\right) - c_2^* \sin\left(\int k_2 ds\right).$$

*and*

$$g(\alpha, B) = c_1 \sin\left(\int k_2 ds\right) - c_2 \cos\left(\int k_2 ds\right)$$

$$g(\alpha, B^*) + g(\alpha^*, B) = \left[c_1 \cos\left(\int k_2 ds\right) + c_2 \sin\left(\int k_2 ds\right)\right] \int (\Delta k_2 + k_2^*) ds$$
$$+ c_1^* \sin\left(\int k_2 ds\right) - c_2^* \cos\left(\int k_2 ds\right).$$

respectively, where $\vec{\alpha}, \vec{N}, \vec{B}$ and $\vec{\alpha}^*, \vec{N}^*, \vec{B}^*$ are real and dual parts of $\tilde{\alpha}, \tilde{N}$ and $\tilde{B}$ respectively. Here, real parts of (10), (14) and (15) are the conditions for a space curve $\vec{\alpha} = \vec{\alpha}(s)$ with Frenet frame $\{\vec{T}, \vec{N}, \vec{B}\}$ and curvatures $k_1$ and $k_2$ to be a normal curve in Euclidean space $E^3$.

Also, we see that $g(\vec{\alpha}, \vec{N}^*) + g(\vec{\alpha}^*, \vec{N}) = \dfrac{k_1^*}{k_1^2}$, i.e., dual part of equation (10) is equal to dual part of equation (14) with opposite sign. So, we can give the following corollary:

**Corollary 3.2.** Let $\tilde{\alpha} = \tilde{\alpha}(s)$ be a dual space curve with curvatures $\overline{k}_1(s) \neq 0$, $\overline{k}_2(s) \neq 0$. Then $\tilde{\alpha}$ is a normal curve in $D^3$ if and only if the following equality holds

$$\frac{1}{\overline{k}_1(s)} = \frac{1}{k_1} - \varepsilon\left(g(\vec{\alpha}, \vec{N}^*) + g(\vec{\alpha}^*, \vec{N})\right). \tag{16}$$

## 4. Dual Spherical Curves in Dual Space $D^3$

In this section, we characterize dual space curves which lie on dual sphere $\tilde{S}^2(\overline{r})$ with radius $\overline{r}$ and center $\overline{c}$.

**Theorem 4.1.** Let $\tilde{\alpha} = \tilde{\alpha}(s)$ be a dual space curve. Then $\tilde{\alpha}$ lies on dual sphere $\tilde{S}^2(\overline{r})$ if and only if

$$\left(\frac{1}{\overline{k}_1}\right)^2 + \left[\left(\frac{1}{\overline{k}_1}\right)' \frac{1}{\overline{k}_2}\right]^2 = \overline{r}^2, \tag{17}$$

holds.

**Proof:** Assume that $\tilde{\alpha}$ lies on dual sphere $\tilde{S}^2(\overline{r})$ with radius $\overline{r}$ and let the center of $\tilde{S}^2(\overline{r})$ be origin 0. Then from the dual position vector we have
$$g(\tilde{\alpha}, \tilde{\alpha}) = \overline{r}^2.$$
The differentiations of this equality with respect to dual arc length $\overline{s}$ give first
$$g(\tilde{\alpha}, \tilde{T}) = 0, \tag{18}$$
and then
$$g(\tilde{\alpha}, \tilde{N}) = -\frac{1}{\overline{k}_1}. \tag{19}$$
and the derivation of last equality gives us
$$g(\tilde{\alpha}, \tilde{B}) = -\left(\frac{1}{\overline{k}_1}\right)' \frac{1}{\overline{k}_2}. \tag{20}$$
Then, from (19) and (20) we can write
$$\tilde{\alpha} = -\frac{1}{\overline{k}_1}\tilde{N} - \left(\frac{1}{\overline{k}_1}\right)' \frac{1}{\overline{k}_2}\tilde{B}. \tag{21}$$

Since we have that the radius of the sphere is $\bar{r}^2 = \|\tilde{\alpha} - 0\|$, we obtain that

$$\left(\frac{1}{\bar{k}_1}\right)^2 + \left[\left(\frac{1}{\bar{k}_1}\right)' \frac{1}{\bar{k}_2}\right]^2 = \bar{r}^2,$$

which completes the proof.

Conversely, assume that regular $C^4$-dual space curve $\tilde{\alpha}$ satisfies the condition of the theorem. Let us consider the parametrized dual curve $\tilde{\alpha} = \tilde{c}(s)$ defined by

$$\tilde{c}(s) = \left(\tilde{\alpha} + \frac{1}{\bar{k}_1} \tilde{N} + \left(\frac{1}{\bar{k}_1}\right)' \frac{1}{\bar{k}_2} \tilde{B}\right), \quad (22)$$

and the function $\bar{r}(s)$ defined by

$$\bar{r}^2 = [\tilde{\alpha} - \tilde{c}]^2 = \left(\frac{1}{\bar{k}_1}\right)^2 + \left[\left(\frac{1}{\bar{k}_1}\right)' \frac{1}{\bar{k}_2}\right]^2, \quad (23)$$

If we differentiate (22) and (23) and make use of dual Frenet formulae, the results are $\tilde{c}' = 0$, $\bar{r}' = 0$, respectively. Therefore, the parametrized curve $\tilde{\alpha} = \tilde{c}(s)$ reduces to a point $\bar{c}$ and the function $\bar{r}(s)$ is a dual constant $\bar{r}$. Hence by (23), $\tilde{\alpha}(s)$ lies on dual sphere $\tilde{S}^2(\bar{r})$ with center $\bar{c}$ and radius $\bar{r}$.

***Corollary 4.1.*** *The real and dual parts of* (17) *are given by*

$$r^2 = \left(\frac{1}{k_1}\right)^2 + \left[\left(\frac{1}{k_1}\right)' \frac{1}{k_2}\right]^2$$

*and*

$$rr^* = -\frac{k_1^*}{k_1^3} - \left[\left(\frac{1}{k_1}\right)'\right]^2 \frac{k_2^*}{k_2^3} + \left(\frac{1}{k_1}\right)' \left(\frac{k_1^*}{k_1^2}\right)' \frac{1}{k_2^2},$$

*respectively. Here* $r^2$, *which is the real part of* (17), *characterizes a space curve* $\vec{\alpha} = \vec{\alpha}(s)$ *with Frenet frame* $\{\vec{T}, \vec{N}, \vec{B}\}$ *and curvatures* $k_1$ *and* $k_2$ *which lies on real sphere* $S^2(r)$ *with radius* $r$ *in* $E^3$ [8].

Equation (21) shows that $\tilde{\alpha}$ is a normal curve in $D^3$. So, we can give the following corollary:

***Corollary 4.2.*** *Let* $\tilde{\alpha} = \tilde{\alpha}(s)$ *be a dual space curve. Then* $\tilde{\alpha}$ *lies on dual sphere* $\tilde{S}^2(\bar{r})$ *if and only if* $\tilde{\alpha}$ *is a normal curve in* $D^3$.

**6. Conclusions**
In this study, the characterizations of dual normal and dual spherical curves are given in dual space $D^3$. Also, it is observed that dual normal curves are also dual spherical curves.